\documentclass[letterpaper, 10 pt, conference]{ieeeconf}  

\IEEEoverridecommandlockouts                              

\usepackage{graphicx} %

\usepackage{amsmath} 
\usepackage{amssymb}  

\usepackage{siunitx}

\title{\LARGE \bf
Online Optimization for Controlling HVAC Systems
}

\author{Tong Wu$^{1,2}$, Dominik Baumann$^{2}$, Angela Fontan$^{1}$, Nicola Bastianello$^{3}$
\thanks{$^{1}$KTH Royal Institute of Technology, Stockholm, Sweden}%
\thanks{$^{2}$Aalto University, Aalto, Finland}%
\thanks{$^{3}$University College London, United Kingdom
        {\tt\small n.bastianello@ucl.ac.uk}}%
}

\begin{document}

\maketitle
\thispagestyle{empty}
\pagestyle{empty}

\begin{abstract}
Heating, ventilation, and air-conditioning (HVAC) systems account for a significant portion of buildings' energy consumption, motivating the design of energy- and cost-efficient control strategies. However, conventional strategies employ simple, reactive control, without incorporating a model. Therefore, in this paper we adopt an economic model predictive control framework and focus on designing online algorithms to solve the control input design problems with a limited computational budget. In particular, we propose two tailored online algorithms based on projected gradient descent, with a custom projection routine, and the projected primal-dual method. We discuss the pros and cons of both algorithm designs, focusing in particular on computational efficiency. We further compare them through numerical experiments simulating the control of an HVAC system over the course of a day. The results show that computationally efficient online solvers for HVAC systems control are indeed possible. Moreover, the performance closely matches that of a baseline algorithm with unlimited computational budget, in terms of reference temperature tracking and comfort temperature range enforcement.
\end{abstract} 

\section{Introduction}
Global energy consumption by end-users exceeded $450 \times 10^{18}~\unit{\joule}$ in 2024, representing an increase of approximately $25 \times 10^{18}~\unit{\joule}$ over 2019 levels \cite{iea2025energyefficiency}.
Buildings account for around $30\%$ of this increasing energy demand \cite{iea2025energyefficiency}.
For example, heating, ventilation, and air-conditioning (HVAC) systems account for approximately $59\%$ of residential buildings' and $36\%$ of commercial buildings' energy consumption in the United States \cite{perez2008review,kong_hvac_2022}.
Improving the operation of HVAC systems is therefore an important means of reducing building energy use and operating costs.

Smart HVAC control strategies could thus reduce energy use and costs, however, conventional strategies rely on simple approaches with low implementation cost, including on/off control, rule-based control, and proportional-integral-derivative (PID) control \cite{maddalena_data_2020}.
However, these methods are primarily reactive, as they do not exploit a model of the system's evolution and predictions of future disturbances, which include fluctuations in outdoor temperature, occupancy , and solar irradiance \cite{afram2014theory,farjadnia_assessing_2026,farjadnia2026generalizability}.
Therefore, model predictive control (MPC) has been proposed as an alternative framework for controlling HVAC systems \cite{drgovna2020all}. The use of MPC allows to integrate predictions of the building's thermal dynamics and disturbances over a finite horizon in order to improve control efficiency, while also allowing to impose comfort and operational constraints.

Beyond operational efficiency and user comfort, energy cost is also an important factor. For this reason, the economic MPC (E-MPC) framework has been further proposed to control HVAC systems \cite{zheng2024economic,boutahri2026economic}.
Indeed, E-MPC allows incorporating the time-varying price of electricity directly in the objective function alongside temperature control goals and constraints. In this paper, we thus adopt the E-MPC framework.

Despite the advantages of the E-MPC framework, in practice it requires solving constrained optimization problems at every sampling instant. The resulting computational burden therefore can become a practical limitation, especially in light of two important factors: thermal inertia and hardware resources.
Buildings, indeed, might be characterized by high thermal inertia, which means that the application of control inputs only very slowly regulates the temperature. It is therefore important to reduce the ratio between computation time and control input application as much as possible, so that temperature can be more effectively regulated during sampling intervals.
The second factor are the computational resources available to implement HVAC systems controllers. These resources are in fact often limited, and computation-heavy control input design might not be possible \cite{parisio2014control}, especially if the computation time exceeds the sampling interval.
Therefore, in this paper we focus on designing efficient solvers for the HVAC E-MPC problem.

From an optimization perspective, HVAC E-MPC consists of a sequence of constrained optimization problems whose parameters change over time in response to fluctuations in energy prices, occupancy, outside temperatures, and past control inputs.
The E-MPC problem therefore can be interpreted naturally as an online optimization problem \cite{dall2020optimization,simonetto2020time}, and more specifically as feedback optimization \cite{colombino2019online}.
Unlike static optimization problems, in online optimization the computational time is limited to the interval between the formulation of consecutive problems. The challenge then is to design algorithms that can approximate the problems' solutions well, but in a limited time.

In this paper, therefore, we focus on the design of online optimization algorithms tailored to the solution of HVAC E-MPC problems with limited computational resources/time.
In particular, we offer the following contributions:
\begin{itemize}
    \item We formulate the E-MPC framework for the control of HVAC systems as an online, constrained optimization problem, using a simple but computationally efficient thermal model of buildings.

    \item We propose two online algorithms tailored to solve this problem: one based on the projected gradient descent, with a customized and efficient projection routine; and one based on the projected primal-dual method, which enforces constraints differently.

    \item We evaluate the two methods with a limited computational budget, both in terms of sub-optimality (that is, how closely they track the sequence of problems' solutions), and in terms of temperature regulation (that is, reference temperature tracking and violation of comfort temperature ranges), providing insights on selecting which algorithm to apply. We also compare them to a baseline algorithm with unlimited computational budget, showing that a limited budget can still achieve strong performance.
\end{itemize}

\section{HVAC System Control as an Online Optimization Problem}

\subsection{Building thermal model}
We model the dominant thermal dynamics of a room within a building with a second-order resistance-capacitance model \cite{afroz2018modeling,yang2024rc}, which allows for physical interpretability and computational tractability.
Representing the indoor air and  building wall as two thermal capacitances, $C_\mathrm{in}$ and $C_\mathrm{w}$ respectively, we can write the energy-balance equations:
\begin{equation}
\begin{aligned}
C_{\mathrm{in}} \dot T_{\mathrm{in}} &=
\frac{T_{\mathrm{w}}-T_{\mathrm{in}}}{R_{\mathrm{iw}}}
+
\frac{T_{\mathrm{out}}-T_{\mathrm{in}}}{R_{\mathrm{iwin}}} \\
&\qquad +
C N_{\mathrm{p}}
+
G A_{\mathrm{win}} I
+
Q_{\mathrm{ctrl}},
\\
C_{\mathrm{w}} \dot T_{\mathrm{w}} &=
\frac{T_{\mathrm{in}}-T_{\mathrm{w}}}{R_{\mathrm{iw}}}
+
\frac{T_{\mathrm{out}}-T_{\mathrm{w}}}{R_{\mathrm{we}}},
\end{aligned}
\label{eq:thermal_balance}
\end{equation}
where $T_{\mathrm{in}}$, $T_{\mathrm{w}}$, and $T_{\mathrm{out}}$ denote the indoor-air temperature, average wall temperature, and outdoor-air temperature ($\unit{\degreeCelsius}$), respectively; $R_{\mathrm{iw}}$, $R_{\mathrm{iwin}}$, and $R_{\mathrm{we}}$ denote the
thermal resistances between the indoor air and wall, indoor air and outdoor environment through the windows, and wall and outdoor environment. The internal heat gain is modeled using the number of occupants $N_{\mathrm{p}}$ and the heat emission per occupant $C$, while $I$, $G$, and $A_{\mathrm{win}}$ denote the solar irradiance, solar heat-gain coefficient, and window area. The control input $Q_{\mathrm{ctrl}}$ represents the aggregate HVAC thermal power supplied to the zone, measured in watts ($\unit{\watt}$).

Defining the state, control input, and disturbance vector as
\begin{equation}
x =
\begin{bmatrix}
T_{\mathrm{in}} \\
T_{\mathrm{w}}
\end{bmatrix},
\qquad
u = Q_{\mathrm{ctrl}},
\qquad
w =
\begin{bmatrix}
T_{\mathrm{out}} \\
N_{\mathrm{p}} \\
I
\end{bmatrix},
\label{eq:thermal_vectors}
\end{equation}
where the disturbances include external temperature, occupancy and solar irradiance, we can then write~\eqref{eq:thermal_balance} as the compact continuous-time dynamics:
$$
\dot x = A x + B u + E w.
$$
Using forward Euler discretization with sampling interval
$\Delta t$, the discrete-time building thermal model is
\begin{equation}
x_{k+1}=A_{\mathrm{d}}x_k+B_{\mathrm{d}}u_k+E_{\mathrm{d}}w_k,
\label{eq:thermal_model}
\end{equation}
where $A_{\mathrm{d}} = \begin{bmatrix} a_1 & a_2 \\ a_3 & a_4 \end{bmatrix}$,
\begin{align*}
    a_1 &= 1-\Delta t\left(
\frac{1}{C_{\mathrm{in}}R_{\mathrm{iw}}}
+\frac{1}{C_{\mathrm{in}}R_{\mathrm{iwin}}}
\right) \qquad
    a_2 = \frac{\Delta t}{C_{\mathrm{in}}R_{\mathrm{iw}}} \\
    a_3 &= \frac{\Delta t}{C_{\mathrm{w}}R_{\mathrm{iw}}} \qquad
    a_4 = 1-\Delta t\left(
\frac{1}{C_{\mathrm{w}}R_{\mathrm{iw}}}
+\frac{1}{C_{\mathrm{w}}R_{\mathrm{we}}}
\right)
\end{align*}
\begin{equation*}
B_{\mathrm{d}} =
\begin{bmatrix}
\frac{\Delta t}{C_{\mathrm{in}}}\\
0
\end{bmatrix},
\qquad
E_{\mathrm{d}} =
\Delta t
\begin{bmatrix}
\frac{1}{C_{\mathrm{in}}R_{\mathrm{iwin}}}
&
\frac{C}{C_{\mathrm{in}}}
&
\frac{GA_{\mathrm{win}}}{C_{\mathrm{in}}}
\\
\frac{1}{C_{\mathrm{w}}R_{\mathrm{we}}}
&
0 & 0
\end{bmatrix}.
\end{equation*}

\subsection{Economic MPC formulation}\label{subsec:empc-formulation}
Having defined the building thermal model~\eqref{eq:thermal_model}, we now formulate the HVAC control problem as an Economic MPC (E-MPC) problem. We select an economic formulation as it allows us to incorporate time-varying electricity prices directly into the control objective \cite{zheng2024economic,boutahri2026economic}, thus adapting it accordingly.

The goal is to control the internal temperature to a reference $T_{\mathrm{ref}}$ temperature, while accounting for the wall temperature $T_{\mathrm{w}}$ as an internal state capturing thermal inertia.
Let $w_k$ and $c_k$ denote the disturbance and electricity-price forecasts at prediction step $k$, respectively, over a prediction horizon of length $H$. These forecasts can be generated using historical data to build a model of both disturbances and price fluctuations.
Then the E-MPC problem at sampling time $k$ is:
\begin{equation}
\begin{aligned}
&\min_{\{x_{h+1},u_h\}_{h=k}^{H+k-1}}
\quad
\sum_{h=k}^{H+k-1}
\Bigg[
\frac{\gamma_{\mathrm{trk}}}{2}
\|x_{h+1}-x_{\mathrm{ref}}\|_{Q}^{2}
\\
& \hspace{3cm}
+\frac{\gamma_\mathrm{ctrl}}{2}\|u_h\|_{R}^{2}
+\gamma_{\mathrm{price}}c_h u_h
\Bigg]
\\
\mathrm{s.t.}\quad
&
x_{h+1}=A_{\mathrm{d}}x_h+B_{\mathrm{d}}u_h+E_{\mathrm{d}}w_h, \\ & \hspace{3cm} h=k,\ldots,H+k-1,
\\
&
T_{\mathrm{in,min}}
\leq T_{\mathrm{in},h}
\leq T_{\mathrm{in,max}},
\quad h=k+1,\ldots,H+k,
\\
&
u_{\min}\leq u_h\leq u_{\max},
\quad h=k,\ldots,H+k-1,
\end{aligned}
\label{eq:empc}
\end{equation}
where $Q$ and $R$ are the state- and input-weighting matrices, and $\gamma_{\mathrm{trk}}$, $\gamma_\mathrm{ctrl}$, and
$\gamma_{\mathrm{price}}$  are the nonnegative weights that determine the relative importance of reference tracking, control effort, and economic cost. We remark that the priority is achieving the internal temperature reference ($T_\mathrm{in} \to T_\mathrm{ref}$), suggesting a choice of $Q$ whose first diagonal term is larger (see Table~\ref{tab:simulation_mpc_configuration}).
The constraints in~\eqref{eq:empc} enforce adherence to the thermal model, the temperature comfort bounds $T_{\mathrm{in,min}}$, $T_{\mathrm{in,max}}$ set by the users, and actuator constraints on the inputs.

We can now write~\eqref{eq:empc} in the compact form of a convex quadratic program. Collecting the predicted states and control inputs over the
horizon in the vector:
\begin{equation}
z_k =
\begin{bmatrix}
x_{k+1}^\top &
\cdots &
x_{H+k}^\top &
u_k &
\cdots &
u_{H+k-1}
\end{bmatrix}^{\top} \in \mathbb{R}^{3H},
\label{eq:stacked_z}
\end{equation}
we can indeed reformulate~\eqref{eq:empc} as:
\begin{equation}
\begin{aligned}
z_k^\star
=
\arg\min_{z}\quad &
\frac{1}{2}z^\top Pz+q_k^\top z
\\
\mathrm{s.t.}\quad &
A_{\mathrm{eq}}z=b_{\mathrm{eq},k},
\\
&
z_{\min}\leq z\leq z_{\max}.
\end{aligned}
\label{eq:time_varying_qp}
\end{equation}
where $P$, $q_k$ collect quadratic and linear terms of the control objective, $A_{\mathrm{eq}}$, $b_{\mathrm{eq},k}$ define the thermal model over the horizon, and $z_{\min}$, $z_{\max}$ collect the lower and upper bounds on the predicted indoor temperature and control input.
Here we assume to have selected positive definite $Q$, $R$ matrices, and
$\gamma_{\mathrm{trk}},\gamma_{\mathrm{ctrl}}>0$. This ensures that the Hessian $P$ is positive definite, and that each QP~\eqref{eq:time_varying_qp} admits a unique minimizer $z_k^\star$.

\subsection{Online optimization formulation}\label{subsec:online-formulation}
As outlined in the previous section, the HVAC E-MPC problem is an online quadratic program, which depends on the state ($x_k$) and disturbance and electricity price forecasts ($w_k$, $c_k$) at sampling time $k$ (with new samples being taken at $\Delta t$ intervals, see~\eqref{eq:thermal_model}).
We note that $P$, $A_{\mathrm{eq}}$, and $z_{\min}$, $z_{\max}$, remain constant over time $k$, as the prediction horizon, thermal model, and temperature and control constraints remain fixed.

\smallskip

Controlling the HVAC system to reach the desired internal temperature range therefore requires solving an online optimization problem \cite{dall2020optimization,simonetto2020time}, that is, of the sequence of problems in~\eqref{eq:time_varying_qp}.
In principle, one might solve each problem in~\eqref{eq:time_varying_qp} in isolation; however, the HVAC use case requires a more computationally efficient strategy. Indeed, the time available to compute the control input $u_k$ to be applied during the interval $[k, k + \Delta t)$ is limited by two factors: \textit{thermal inertia} and \textit{computational resources}.
First of all, buildings can have high thermal inertia, meaning that each input $u_k$ should be applied for as close to $\Delta t$ as possible in order to successfully control the system to the desired temperature reference.
Secondly, HVAC controllers might operate on low-cost hardware with limited computational resources, excluding solvers for~\eqref{eq:time_varying_qp} that need large computational budgets, as well as more complicated problem formulations (that is, more complicated thermal models than~\eqref{eq:thermal_balance}).
The goal of this paper then is to design efficient online algorithms that can solve~\eqref{eq:time_varying_qp} with a limited computational budget.

\section{Proposed Algorithms}\label{sec:algorithms}
In this section, we investigate two online convex optimization methods for solving~\eqref{eq:time_varying_qp} under limited
computational resources. At each MPC sampling instant, updated measurements and forecasts define a new optimization problem, whose solution determines the next HVAC control input. To provide this input in a timely manner, the algorithms perform a limited number of optimization iterations rather than requiring full convergence at each sampling instant.
Denoting the objective and constraint set of~\eqref{eq:time_varying_qp} by
$
    f_k(z)=\frac{1}{2}z^\top Pz+q_k^\top z
$
and
$
    \mathcal{C}_k =
    \left\{
    z \,\middle|\,
    A_{\mathrm{eq}}z=b_{\mathrm{eq},k},\;
    z_{\min}\leq z\leq z_{\max}
    \right\},
$
we can rewrite the problem as $z_k^\star = \arg\min_z f_k(z) \ \text{s.t.} \ z \in \mathcal{C}_k$.
The next two sections describe the proposed iterative algorithms for solving this problem.

\subsection{Algorithm 1: Projected gradient method}
The first algorithm we explore is the \textit{online projected gradient descent} (PG) \cite{dall2020optimization}. Let $N_{\mathrm{iter}} \in \mathbb{N}$ be the computational budget, that is, the number of iterations of PG that we apply, and $\ell$ the iteration counter; then the algorithm is characterized by the updates:
\begin{subequations}\label{eq:pg}
\begin{align}
    \tilde{z}_k^{\ell+1} &= z_k^\ell - \alpha_{\mathrm{PG}} \left(Pz_k^\ell+q_k\right) \label{eq:pg_gradient} \\
    z_k^{\ell+1} &= \Pi_{\mathcal{C}_k}(\tilde{z}_k^{\ell+1}), \quad \ell = 0, \ldots, N_{\mathrm{iter}}-1,
\end{align}
\end{subequations}
where $\alpha_{\mathrm{PG}}>0$ denotes the step-size, and the projection onto the convex set $\mathcal{C}_k$ is~\cite{parikh2014proximal}
\begin{equation}
    \Pi_{\mathcal{C}_k}(\tilde z)
    =
    \arg\min_{z\in\mathcal{C}_k}
    \frac{1}{2}\|z-\tilde z\|_2^2.
    \label{eq:euclidean_projection}
\end{equation}

In the HVAC control problem~\eqref{eq:time_varying_qp}, the constraint set is the intersection of two sets
$
    \mathcal{C}_k
    =
    \mathcal{C}_{\mathrm{box}}
    \cap
    \mathcal{C}_{\mathrm{eq},k},
$
with
\begin{align*}
    \mathcal{C}_{\mathrm{box}}
    &=
    \left\{
    z \,\middle|\,
    z_{\min}\leq z\leq z_{\max}
    \right\},\\
    \mathcal{C}_{\mathrm{eq},k}
    &=
    \left\{
    z \,\middle|\,
    A_{\mathrm{eq}}z=b_{\mathrm{eq},k}
    \right\},
\end{align*}
and thus the projection does not have a closed form. Instead of directly solving the quadratic program~\eqref{eq:euclidean_projection}, we propose to approximately compute the projection using the \textit{alternating projection} method \cite{reich_projection_2015}. As the name suggests, the idea is to alternate between projections onto the two sets; initializing $w_k^0=\tilde z_k^{\ell+1}$, we do:
\begin{equation}
\begin{aligned}
    \bar w_k^{j+1}
    &=
    \Pi_{\mathcal{C}_{\mathrm{eq},k}}(w_k^j),\\
    w_k^{j+1}
    &=
    \Pi_{\mathcal{C}_{\mathrm{box}}}(\bar w_k^{j+1}),
    \qquad
    j=0,\ldots,N_{\mathrm{proj}}-1,
\end{aligned}
\label{eq:alternating_projection}
\end{equation}
where $N_{\mathrm{proj}}$ denotes the number of iterations.
The advantage of this approach is that the two projections have a closed form, which is computationally efficient; indeed:
\begin{equation}
    \Pi_{\mathcal{C}_{\mathrm{box}}}(z)
    =
    \min\!\left(
    \max(z,z_{\min}),z_{\max}
    \right),
    \label{eq:box_projection}
\end{equation}
and
\begin{equation}
\begin{split}
    \Pi_{\mathcal{C}_{\mathrm{eq},k}}(\tilde z)
    =
    \tilde z
    -
    A_{\mathrm{eq}}^\top
    \left(
    A_{\mathrm{eq}}A_{\mathrm{eq}}^\top
    \right)^{-1}
    \left(
    A_{\mathrm{eq}}\tilde z-b_{\mathrm{eq},k}
    \right),
\end{split}
\label{eq:eq_projection}
\end{equation}
since $A_{\mathrm{eq}}$ has full row rank. Notice that the inverse $\left( A_{\mathrm{eq}}A_{\mathrm{eq}}^\top \right)^{-1}$ can be computed once offline, thus requiring only matrix-vector multiplications for~\eqref{eq:eq_projection}.

The alternating projection method is guaranteed to yield the projection of $z_k^{\ell+1}$ onto the intersection asymptotically, so $N_{\mathrm{proj}}$ iterations only yield an approximation thereof. Nonetheless, in practice this approximation is very accurate even with small $N_{\mathrm{proj}}$, given the relatively simple structure of the constraint set, as shown in section~\ref{sec:numerical}.

\subsection{Algorithm 2: Projected primal-dual method}
The PG method discussed in the previous section requires the projection onto the intersection of two sets. An alternative to that is to apply the projected primal-dual (PPD) method~\cite{bernstein2019online}: the box constraints are enforced directly through projection, while the equality constraints are enforced through a penalty term.
Let $\lambda_k^\ell \in \mathbb{R}^{2H}$ denote the dual variable associated to the equality constraints; then the algorithm is characterized by the updates:
\begin{subequations}
\begin{align}
    z_k^{\ell+1}
    &=
    \Pi_{\mathcal{C}_{\mathrm{box}}}
    \Big[
    z_k^\ell
    -
    \alpha_{\mathrm{PPD}}
    \big(
    Pz_k^\ell+q_k
    +A_{\mathrm{eq}}^\top\lambda_k^\ell \label{eq:ppd_primal}\\
    &\hspace{1.3cm} +\rho A_{\mathrm{eq}}^\top
    (A_{\mathrm{eq}}z_k^\ell-b_{\mathrm{eq},k})
    \big)
    \Big], \nonumber\\
    \lambda_k^{\ell+1}
    &=
    \lambda_k^\ell
    +
    \beta_{\mathrm{PPD}}
    \left(
    A_{\mathrm{eq}}z_k^{\ell+1}
    -
    b_{\mathrm{eq},k}
    \right),
    \label{eq:ppd_dual}
\end{align}
\end{subequations}
where $\alpha_{\mathrm{PPD}}, \ \beta_{\mathrm{PPD}}>0$ are the primal and dual step-sizes, and $\rho>0$ the equality constraints penalty weight.
This algorithm prioritizes enforcement of the box constraints (temperature comfort bounds) at every iteration, to adherence to the thermal dynamics, since the equality constraints are only satisfied asymptotically.

\smallskip

To summarize, PG and PPD employ different mechanisms for handling the constraints in the E-MPC problem: PG approximates the projection onto the intersection of box and equality constraints, while PPD projects only onto the box constraints and iteratively penalizes the equality constraints violation.
Their different designs result in different trade-offs between optimization accuracy, constraints violation, and computational complexity. We next evaluate these trade-offs numerically.

\section{Numerical Results}\label{sec:numerical}
In this section, we compare the performance of the algorithms presented in section~\ref{sec:algorithms} against each other and a reference solver (based on OSQP \cite{stellato2020osqp}).

\subsection{Simulation setup and evaluation metrics}
We simulate the control of an HVAC system over a $24$-hour period, with a sampling interval of $\Delta t = 15~\unit{\minute}$.
We use synthetic disturbances (external temperature, occupancy, solar irradiance) and price profiles for~\eqref{eq:time_varying_qp}. In particular, these model a typical residential building profile, with residents being present mostly in the morning and evening, and out during most of the day.
The main simulation and solver parameters are summarized in Table~\ref{tab:simulation_mpc_configuration}.
\begin{table}[!ht]
\centering
\caption{Main simulation and solver parameters.}
\label{tab:simulation_mpc_configuration}
\small
\begin{tabular}{lll}
\hline
Parameter & Symbol & Value \\
\hline
Sampling time
& $\Delta t$
& $15~\unit{\minute}$ \\
Prediction horizon
& $H$
& $20$ steps \\
Simulation length
& $N_{\mathrm{sim}}$
& $96$ steps \\
Initial state
& $x_0$
& $[21,\;20]^\top\unit{\degreeCelsius}$ \\
Reference temp.
& $T_{\mathrm{ref}}$
& $22\unit{\degreeCelsius}$ \\
State weight
& $Q$
& $\mathrm{diag}(5000,1)$ \\
Input weight
& $R$
& $10^{-5}$ \\
Objective weights
& $\gamma_{\mathrm{trk}},\gamma_\mathrm{ctrl},\gamma_{\mathrm{price}}$
& $1,1,1$ \\
Comfort bounds
& $T_{\mathrm{in,min}},T_{\mathrm{in,max}}$
& $20,24\unit{\degreeCelsius}$ \\
HVAC bounds
& $u_{\min},u_{\max}$
& $0,5000~\unit{\watt}$ \\
PG step-size
& $\alpha_{\mathrm{PG}}$
& $3\times10^{-6}$ \\
PPD step-sizes
& $\alpha_{\mathrm{PPD}},\beta_{\mathrm{PPD}}$
& $5\times10^{-5},3\times10^{-5}$ \\
PPD parameter
& $\rho$
& $1$ \\
Computation time
& $t_{\mathrm{budget}}$
& $2~\unit{\milli\second}$ \\
\hline
\end{tabular}
\end{table}
We remark that the algorithms' step-sizes were selected by manual tuning, and both PG and PPD are given a very tight computational budget of $2~\unit{\milli\second}$ at each time $k$ (as motivated in section~\ref{subsec:online-formulation}). The horizon characterizing the E-MPC problem is selected at $H = 20$.%
\footnote{We also remark that, as discussed in section~\ref{subsec:empc-formulation}, the most important objective is for the internal temperature to reach $T_{\mathrm{ref}}$. This motivates the choice of a much larger weight in the first diagonal element of $Q$.}

We evaluate the performance of the proposed algorithms in two ways.
The first is by evaluating the \textit{instantaneous sub-optimality} defined by
\begin{equation}
    e_k =
    \left\|
    z_k -
    z_k^\star
    \right\|_2 ,
    \label{eq:solver_tracking_error}
\end{equation}
where $z_k$ is the output of the algorithm, and $z_k^\star$ the solution to~\eqref{eq:time_varying_qp}. Notice that $z_k^\star$ depends on the previous sequence of control inputs; hence, the cost, constraints, and the optimal solution at time $k$ are different between PG and PPD (or between instances of the same algorithm with different step-sizes). Concretely, $z_k^\star$ is computed using OSQP \cite{stellato2020osqp}.

The second performance metric is a comparison with an alternative solver that has unlimited computational budget. This solver can solve~\eqref{eq:time_varying_qp} exactly, and thus selects the optimal control input at every time $k$. As discussed above, the sequence of $z_k^\star$'s for each solver is different; therefore, we compare against the unlimited budget solver by evaluating the indoor temperature mean absolute error:
\begin{equation}
    \mathrm{MAE} =
    \frac{1}{N_{\mathrm{sim}}}
    \sum_{k=0}^{N_{\mathrm{sim}}-1}
    \left|T_{\mathrm{in},k}-T_{\mathrm{ref}}\right|.
    \label{eq:temperature_mae}
\end{equation}

\subsection{Instantaneous sub-optimality comparison}
\begin{figure}[!ht]
\centering
\includegraphics[width=\linewidth]{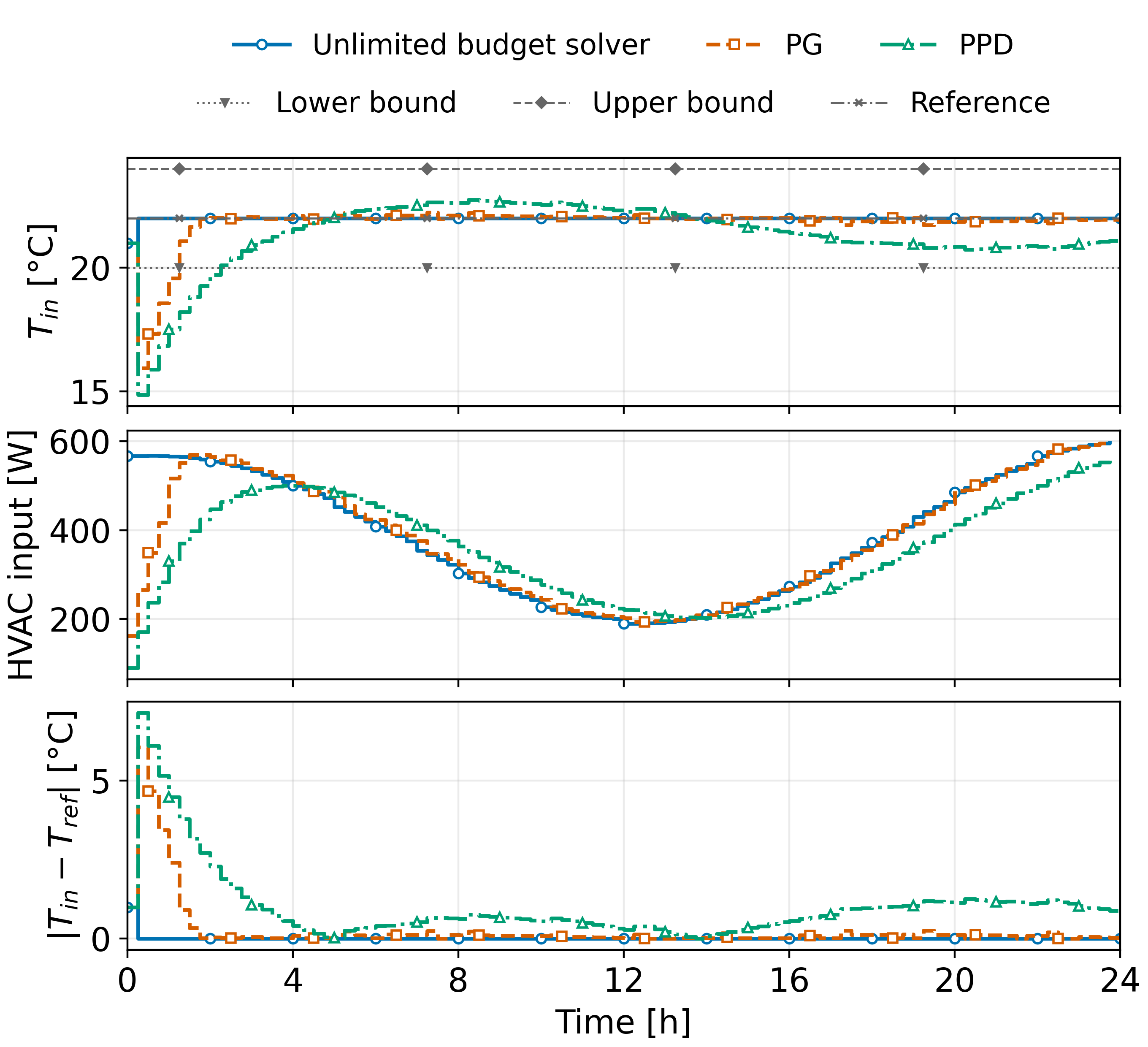}
\caption{Closed-loop HVAC performance of PG and PPD compared with the unlimited budget solver.}
\label{fig:solvers_performance}
\end{figure}

We start by comparing the two proposed algorithms, PG and PPD, in terms of the instantaneous sub-optimality defined in~\eqref{eq:solver_tracking_error}. As mentioned before, both algorithms are given the same computational budget of $2~\unit{\milli\second}$.
Table~\ref{tab:overall_solver_comparison} shows the mean and final sub-optimality $e_k$ achieved by the two algorithms. We also report the mean number of algorithm iterations that can be run during the $2~\unit{\milli\second}$ computational budget.
\begin{table}[!ht]
\centering
\caption{Comparison of instantaneous sub-optimality performance.}
\label{tab:overall_solver_comparison}
\small
\begin{tabular}{lrrr}
\hline
Algorithm
& Mean $e_k$
& Final $e_k$
& Mean iter. \# \\
\hline
PG
& $67.28$
& $2.39$
& $1.45$ \\
PPD
& $670.09$
& $758.54$
& $36.21$ \\
\hline
\end{tabular}
\end{table}
As we can see, PG achieves significantly lower sub-optimality than PPD, even though its iterations are more computationally expensive and thus fewer of them are applied. An iteration of PG is indeed more time-consuming since it integrates the alternating projections to approximate the projection onto the constraint set $\mathcal{C}_k$. PPD, on the other hand, integrates only a simple projection onto the box constraints, while the equality constraints are enforced through a penalty term. As a consequence, PPD is less computationally expensive, but at the cost that it would require a much larger number of iterations to reach similar sub-optimality levels as PG.
Figure~\ref{fig:solvers_performance} further depicts the sub-optimality of the two algorithms over the whole simulation.
\begin{figure}[!ht]
\centering
\includegraphics[width=\linewidth]{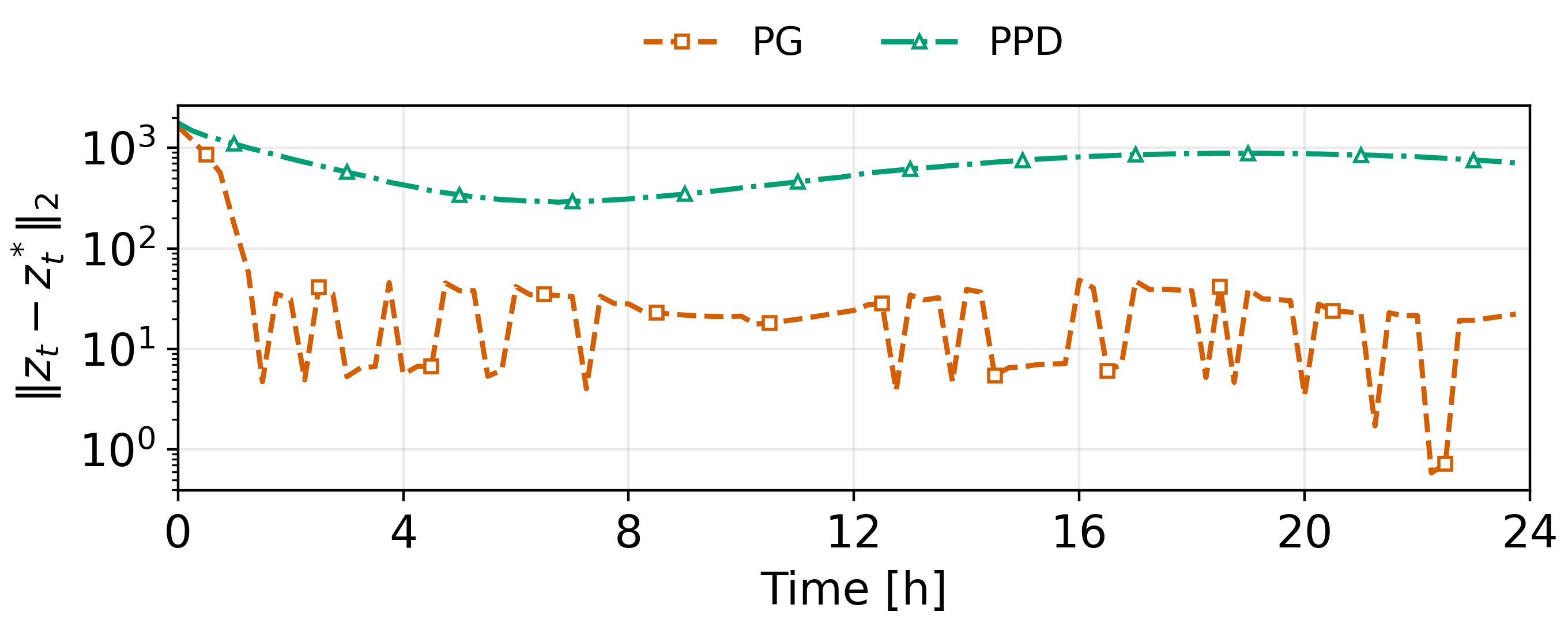}
\caption{Comparison of instantaneous sub-optimality over the simulation.}
\label{fig:solvers_optimality}
\end{figure}
We can see that, after the initial transients, PG stabilizes around $e_k$ values that are two orders of magnitude of what PPD achieves.

These results highlight that how constraints are enforced has a significant impact on the sub-optimality achieved by the algorithms. Therefore, we evaluate now the tracking of the temperature reference $T_{\mathrm{ref}}$ (via the mean absolute error defined in~\eqref{eq:temperature_mae}) alongside the violation of the temperature range constraints ($T_{\mathrm{in,min}}$, $T_{\mathrm{in,max}}$).
The results are reported in Figure~\ref{fig:solvers_performance} and Table~\ref{tab:closed_loop_performance}, compared against the solver with unlimited computational budget as a baseline.
\begin{table}[!ht]
\centering
\caption{Performance in terms of reference temperature and comfort temperature range.}
\label{tab:closed_loop_performance}
\small
\begin{tabular}{lrr}
\hline
Algorithm
& MAE [$\unit{\degreeCelsius}$]
& Max. viol. [$\unit{\degreeCelsius}$] \\
\hline
Unlim. budget sol.
& $0.0105$
& $0$ \\
PG
& $0.2669$
& $4.0663$ \\
PPD
& $1.1078$
& $5.1153$ \\
\hline
\end{tabular}
\end{table}
We can see that the baseline solution achieves zero comfort temperature range violation across the whole simulation (as it exactly satisfies the constraints), while the iterative algorithms PG and PPD exhibit some constraint violation during the initial transient. After this transient, however, both algorithms manage to keep the temperature within the comfort range, with PG tracking the reference temperature more closely than PPD.
We note that the initial transient is caused by initializing the building at $15\unit{\degreeCelsius}$, significantly below $T_{\mathrm{in,min}}$.

The results of this comparison show that the choice of solution algorithm, and particularly how it is designed to enforce the constraints, plays a significant role in tracking the reference temperature and keeping it within the comfort range. In particular, PG, integrating the alternating projections, outperforms PPD.
In the following section, we thus focus on PG alone and evaluate how its performance changes based on the computational budget that is available to it.

\subsection{Effect of computational budget}
As discussed in section~\ref{subsec:online-formulation}, the computational budget available to the algorithm computing the control input depends on two main factors: the thermal inertia and the hardware resources. In this section, we evaluate the impact of changes in the computational budget motivated by \textit{e.g.} the use of state-of-the-art HVAC systems that regulate the temperature more efficiently (and hence, rapidly), or the use of better computational resources.

In particular, we fix the computational budget of PG by setting the number of iterations that can be applied at every time $k$ to compute the control input, selecting:
\begin{equation}
    N_{\mathrm{iter}}\in\{1,10,100\},
\end{equation}
which roughly correspond to $1.4~\unit{\milli\second}$, $14~\unit{\milli\second}$ and $140~\unit{\milli\second}$ of compute time.
All other simulation settings and PG parameters, reported in Table~\ref{tab:simulation_mpc_configuration}, are unchanged.

\begin{figure}[!ht]
\centering
\includegraphics[width=\linewidth]
{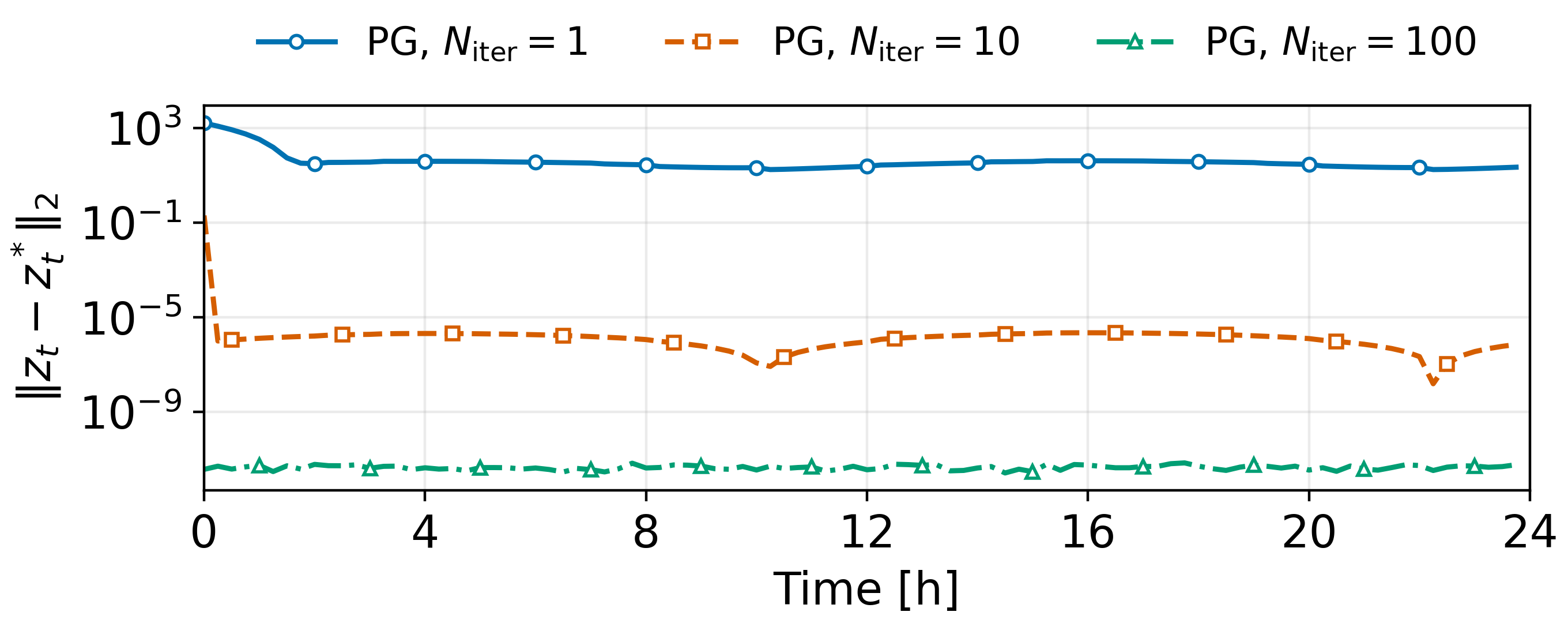}
\caption{Sub-optimality comparison for PG with different computational budgets.}
\label{fig:pg_iterationcompare_optimality}
\end{figure}

Figure~\ref{fig:pg_iterationcompare_optimality} depicts the sub-optimality of PG with different computational budgets across the $24\unit{\hour}$ simulation.
Clearly, increasing the budget reduces sub-optimality, with the mean $e_k$ corresponding to the three budgets being $78.51$, $2.11\times10^{-3}$, and $4.64\times10^{-12}$.

\begin{figure}[!ht]
\centering
\includegraphics[width=\linewidth]
{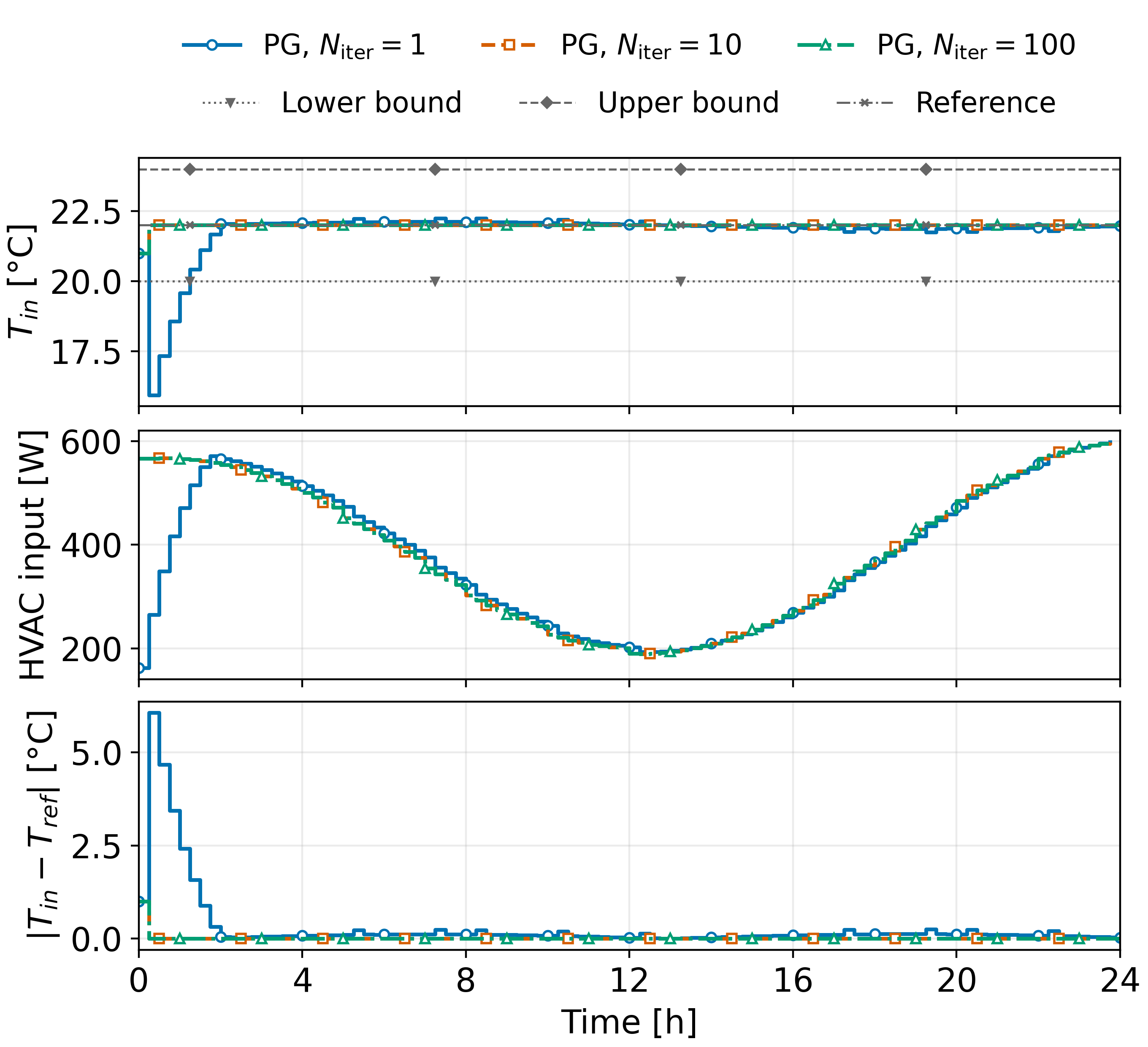}
\caption{Closed-loop performance of PG for different numbers of optimization
iterations.}
\label{fig:pg_iterationcompare_performance}
\end{figure}

Figure~\ref{fig:pg_iterationcompare_performance} depicts the performance also in terms of reference temperature tracking and comfort temperature range violation.
We can see that with one iteration, the
temperature MAE is $0.302\unit{\degreeCelsius}$, and with an initial transient during which the comfort range is violated.
With 10 iterations, the MAE decreases
to $0.01048\unit{\degreeCelsius}$ and no comfort violations occur across the whole simulation.
Further increasing the computational budget, while reducing the sub-optimality, does not, however, significantly reduce the MAE, since the maximum difference between PG-10 and PG-100
temperature trajectories is only
$6.13\times10^{-4}\unit{\degreeCelsius}$.
This indicates that, while a larger computational budget can be useful, this usefulness plateaus after a while. This is an important observation, as it ensures that satisfactory performance in terms of temperature tracking can still be achieved even with limited time (that is, high thermal inertia or poor computational resources).

\smallskip

To summarize the results of these numerical results, enforcing the constraints of~\eqref{eq:time_varying_qp} more strictly, like PG does, has a significant impact on the performance of the HVAC controller, especially in terms of temperature range violation.
Additionally, a modest computational budget can still yield very good performance, satisfying the practical constraints of buildings fitted with HVAC systems.

\section{Conclusions}
In this paper, we focused on the control of heating, ventilation, and air-conditioning (HVAC) systems, with the goal of reducing their energy consumption while guaranteeing user comfort.
To this end, we have adopted the E-MPC framework and discussed its interpretation as an online, constrained optimization problem.
We have proposed two tailored online algorithms to solve this problem, based on projected gradient descent with a custom projection routine and the projected primal-dual method, discussing the pros and cons of both.
We ran numerical experiments to thoroughly compare the two algorithms, both against each other and against a baseline algorithm with unlimited computational budget.
The results showed that computational efficiency does not necessarily come at the cost of worse performance, especially for the projected gradient descent method. This algorithm indeed closely matched the performance of the baseline algorithm in terms of reference temperature tracking and comfort temperature range enforcement.
Future work will focus on improving the computational efficiency of the algorithms, integrating model and thermal disturbances uncertainty, including the effects of occupants' presence and behavior \cite{farjadnia_assessing_2026,farjadnia2026generalizability}.

\addtolength{\textheight}{-12cm}   




\bibliographystyle{IEEEtran}
\bibliography{references}

\end{document}